\documentclass[12pt]{article}
\usepackage{amsmath,amssymb,amsthm,eucal}
\usepackage{hyperref}
\usepackage{color}
\usepackage[normalem]{ulem}
\usepackage{enumitem}

\font\tenrm=cmr10

\font\cmssl=cmss10 at 12 pt

\font\bigss=cmssdc10 scaled 2300

\font\cmsslll=cmss10 at 14 pt

\renewcommand{\b}{\beta}

\renewcommand{\d}{\delta}
\newcommand{\e}{\epsilon}

\newcommand{\g}{\gamma}

\renewcommand{\l}{\lambda}

\newcommand{\G}{\Gamma}

\newcommand{\bR}{\mathbb{R}}

\newcommand\GL{\mathrm{GL}}

\newcommand\SU{\mathrm{SU}}

\renewcommand{\square}{\kern1pt\vbox
               {\hrule height 0.6pt\hbox{\vrule width 0.6pt\hskip 3pt
    \vbox{\vskip 6pt}\hskip 3pt\vrule width 0.6pt}\hrule height0.6pt}
    \kern1pt}

\newcommand{\ra}{\rightarrow}

\DeclareMathOperator\tr{tr\;}

\newcommand{\ol}{\overline}

\newtheorem{Pb}{Problem}

\newtheorem{Th}{Theorem}[section]
\newtheorem{Prop}[Th]{Proposition}

\newtheorem{Cor}[Th]{Corollary}
\newtheorem{Lem}[Th]{Lemma}
\newtheorem{Def}[Th]{Definition}

\newtheorem{open}[Th]{Open problem}
\newcommand{\bP}{\begin{Pb}\ \ }
\newcommand{\eP}{\end{Pb}}
\newcommand{\bt}{\begin{Th}\ \ }
\newcommand{\et}{\end{Th}}
\newcommand{\bp}{\begin{Prop}\ \ }
\newcommand{\ep}{\end{Prop}}
\newcommand{\bc}{\begin{Cor}\ \ }
\newcommand{\ec}{\end{Cor}}
\newcommand{\bl}{\begin{Lem}\ \ }
\newcommand{\el}{\end{Lem}}
\newcommand{\bd}{\begin{Def}\ \ }
\newcommand{\ed}{\end{Def}}

\newcommand{\pf}{\begin{proof}[{\it Proof:\ \ }]}
\newcommand{\epf}{\end{proof}}
\newcommand{\n}{\nabla}

\newcommand{\ot}{\otimes}

\newcommand{\be}{\begin{equation}}
\newcommand{\ee}{\end{equation}}

\newcommand\re[1]{(\ref{#1})}
\newcommand{\arr}{\begin{array}{rlll}}
\newcommand{\ea}{\end{array}}
\newcommand{\bea}{\begin{eqnarray}}
\newcommand{\eea}{\end{eqnarray}}
\newcommand{\bean}{\begin{eqnarray*}}
\newcommand{\eean}{\end{eqnarray*}}

\catcode`@=11
\@addtoreset{equation}{section}
\catcode`@=12

\begin{document}
 \rightline{}
\rightline{06.06.2016}
\vskip 1.5 true cm
\begin{center}
{\bigss  Completeness of hyperbolic centroaffine\\[.5em]
 hypersurfaces}\\[.5em]
\vskip 1.0 true cm
{\cmsslll  V.\ Cort\'es$^1$, M.\ Nardmann$^{1,2}$ and S.\ Suhr$^{1,3}$
} \\[3pt]
$^1${\tenrm
Fachbereich Mathematik\\
und Zentrum f\"ur Mathematische Physik\\
Universit\"at  Hamburg\\
Bundesstra{\ss}e 55\\
D-20146 Hamburg, Germany\\
vicente.cortes@uni-hamburg.de
}
 \\[3pt]
$^2${\tenrm
Fakult\"at f\"ur Mathematik\\
TU Dortmund\\
Vogelpothsweg 87\\
D-44227 Dortmund, Germany\\
marc.nardmann@udo.edu
}
 \\[3pt]
$^3${\tenrm
D\'epartement de math\'ematiques et applications\\
ENS Paris and Universit\'e Paris Dauphine\\
45 rue d'Ulm\\
F-75230 Paris Cedex 05, France \\
suhr@dma.ens.fr
}
\end{center}
\vskip 1.0 true cm
\baselineskip=18pt
\begin{abstract}
\noindent
This paper is concerned with the completeness (with respect to the centroaffine metric)
of hyperbolic centroaffine hypersurfaces
which are closed in the ambient vector space. We show that completeness
holds under generic regularity conditions on the boundary of the convex cone generated by the
hypersurface. The main result is that completeness holds  for hyperbolic components of level sets of
homogeneous cubic polynomials. This implies that every such component defines
a complete quaternionic K\"ahler manifold of negative scalar curvature.
\\[.5em]
\end{abstract}
\tableofcontents

\section*{Introduction}

By a celebrated theorem of Cheng and Yau \cite{CY} a locally strictly convex affine hypersphere
which is closed in the ambient vector space is complete with respect to the Blaschke metric.
Proper affine hyperspheres are precisely the centroaffine hypersurfaces for which the
Blaschke metric coincides with the centroaffine metric (up to a constant factor).
In this paper we investigate the completeness of locally strictly convex centroaffine
hypersurfaces with respect to the centroaffine metric.

Our main motivation
stems from the scalar geometry of $5$-dimensional supergravity as described in \cite{GST}.
The manifolds carrying this geometry are called {\it projective special real manifolds}, see Definition \ref{PSRDef}.
They form a class of hyperbolic (and thus locally strictly convex) centroaffine hypersurfaces.
In Theorem \ref{intrinsicThm} and Definition \ref{intrinsicDef} we give an intrinsic characterization in terms of the underlying
centroaffine geometry. The crucial ingredient is the differential equation \re{fundEqu}  expressing the covariant
derivative of the cubic form in terms of the metric.

Using constructions from supergravity, known as the r-map and the c-map,  it was shown in \cite{CHM},  where these constructions are explicitly described, that every \emph{complete} projective special real
manifold of dimension $n$ gives rise to
a complete quaternionic K\"ahler manifold of negative scalar curvature of dimension $4n+8$.
More specifically, it was shown that  given a complete projective special real manifold of dimension $n$,
the r-map associates with it a complete \emph{projective special K\"ahler domain} (see \cite[p.\ 198]{CHM} for a definition) of real dimension $2n+2$ and that the c-map associates a
complete quaternionic K\"ahler manifold of real dimension $4n+8$ of negative scalar curvature with the latter.
This method was used in \cite{CHM,CDL} to construct new explicit examples of complete quaternionic K\"ahler manifolds of dimension $12$ and $16$. Moreover, a classification of all complete projective special real manifolds
of dimension less than or equal to $2$ was given.
Based on these results, it was observed  \cite[Cor.\ 1]{CDL} that
a projective special real manifold
of dimension less than or equal to $2$ is complete if and only if it is closed and it was asked whether this
property extends to higher dimensions. Here we prove that this is indeed the case, see Theorem \ref{mainThm}. This gives
a powerful method for the verification of the completeness of projective special real
manifolds and the corresponding quaternionic K\"ahler manifolds, cf.\ Theorem \ref{main_applThm}.

Let us now summarize the structure of the paper and mention some further results.
In Section \ref{CentroaffSec} we discuss centroaffine structures and centroaffine hypersurfaces from an intrinsic as well as extrinsic point of view.
Our main focus is on locally strictly convex centroaffine hypersurfaces and the relation between:
\begin{enumerate}
\item closedness (the property of being closed in the ambient space),
\item Euclidean completeness (completeness with respect to the metric induced by a Euclidean scalar product in the ambient space) and
\item completeness (with respect to the centroaffine metric).
\end{enumerate}
Section \ref{closeSec} contains some basic results relating these properties. In particular, under natural
assumptions, completeness implies closedness and the latter is equivalent to Euclidean completeness, see Proposition \ref{closedProp}.

After these preliminaries, we
concentrate on Euclidean complete hyperbolic centroaffine hypersurfaces $\mathcal{H}\subset \bR^{n+1}$
in Section \ref{auxSec}. We show that $U = \bR^{>0}\cdot \mathcal{H}$ is an open convex cone, which
is intersected in a relatively compact domain $B=U\cap E\subset E$ by any affine hyperplane $E$ tangent to $\mathcal{H}$.
We equip $U$ with a smooth homogeneous
function $h : U \ra \bR$ of degree $k>1$ such that $\mathcal{H}= \{ p\in U \mathbin{|} h(p)=1\}$ and with a Lorentzian metric
$g_L$ which is a multiple of the Hessian of $h$. We observe that the completeness of $\mathcal{H}$ is equivalent to the
global hyperbolicity of $(U,g_L)$. Then we prove that
$\mathcal{H}$ is complete if there exists $\e\in (0,k)$ such that $f=\sqrt[k-\varepsilon]{h}\big|_B$ is concave, see Lemma
\ref{concaveLemma}. This allows us to prove the completeness if $h$ is a cubic polynomial and
is the key lemma for the proof of Theorem \ref{mainThm} about projective special real manifolds.
As discussed in the last section of the paper, the result for cubic polynomials can not be extended to real analytic functions but might hold for polynomials of higher degree.

In Section \ref{boundarySec} we prove that  a Euclidean complete hyperbolic centroaffine hypersurface
is complete if  the boundary of $U$ satisfies certain regularity assumptions, see Theorem \ref{regThm}.
Furthermore, these conditions are generically satisfied by Theorem \ref{genThm}.

In Section \ref{2ndSec} we specialize to the case of projective special real manifolds. The main
results are the intrinsic characterization developed in Section \ref{intrinsicSec} and the equivalence
of closedness and completeness proven in Section \ref{mainSec}, with the application to
quaternionic K\"ahler manifolds in Section \ref{ApplSec}.

\subsubsection*{Acknowledgements}

This work was partly supported by the German Science Foundation (DFG) under the
Collaborative Research Center (SFB) 676 Particles, Strings and the Early Universe,
as well as by the European Research Council under the European Union's Seventh
Framework Programme (FP/2007-2013) / ERC Grant Agreement 307062.\\
V.C.\ would like to thank Thomas Mohaupt  for collaboration
on the intrinsic geometry of projective special real manifolds
in the context of a related project.  He would also like to thank
Antonio Mart\'{\i}nez and Miguel S\'anchez for discussions
during visits of the University of Granada.\\
Further the authors would like to thank A. Rod Gover for pointing our attention towards
the results in \cite{melrose} and explaining the results of \cite{CG} and \cite{CG2}.

\section{Centroaffine structures} \label{CentroaffSec}

\subsection{Centroaffine hypersurfaces and centroaffine structures}

In this subsection we review some basic notions
from affine differential geometry, see \cite{NS} for
a more detailed discussion.
Let us consider $\bR^{n+1}$ endowed with its canonical
flat connection $\tilde\nabla$ and the parallel volume form
$\det$.

\bd
A hypersurface immersion $\varphi : M \ra \bR^{n+1}$ is called
{\cmssl centroaffine} if the position vector field
$\xi : M \ra \bR^{n+1}$, $p\mapsto \xi_p := \varphi (p)$, is
transversal, that is
for all $p\in M$ we have $\varphi (p)\not\in d\varphi T_pM$.
\ed

A centroaffine  hypersurface immersion $\varphi : M \ra \bR^{n+1}$
induces a torsion-free connection $\n$ and a symmetric tensor field
$g$ on $M$ such that the {\cmssl Gau{\ss} equation}
\be \label{GE}\tilde\n_Xd\varphi Y = d\varphi \n_XY + g(X,Y)\xi\ee
holds for all $X,Y\in \mathfrak{X}(M)$.
Furthermore, we have an induced volume form
$\nu :=\det (\xi , \cdots )$, which is $\n$-parallel, as a consequence of
the {\cmssl Weingarten equation} $\tilde\n_X\xi = d\varphi X$. In these
formulas, $\tilde\n$ denotes the connection
in $\varphi^*T\bR^{n+1}$ induced by the connection $\tilde\n$
in the vector bundle $T\bR^{n+1}$. For simplicity of notation, we
will usually identify
$TM$ with the subbundle $d\varphi TM \subset \varphi^*T\bR^{n+1}$
and drop the isomorphism $d\varphi: TM \ra d\varphi TM$
in the equations of Gau{\ss} and Weingarten.

\bd
The above geometric data $(\n,g,\nu)$ will be called
the {\cmssl induced (centroaffine) data} of the centroaffine hypersurface
immersion $\varphi$. The hypersurface (immersion) is called
{\cmssl nondegenerate} if $g$ is nondegenerate and {\cmssl definite} if $g$ is definite.
More specifically, it is called {\cmssl elliptic}
if $g<0$ and {\cmssl hyperbolic} if $g>0$.
\ed

\noindent
{\bf Remark}: The above definition is consistent with the usual notions of ellipticity and hyperbolicity in
affine differential geometry. In fact, the tensor field $(\tr S)g$ is positive definite in the elliptic case
and negative definite in the hyperbolic case, where $S=-\mathrm{Id}$ is the shape tensor associated with
our choice of transversal vector field $\xi$. This ensures for instance that ellipsoids around $0$ are elliptic and standard hyperboloids are hyperbolic.

\noindent
{\bf Example:}
Let $U\subset \bR^{n+1}$ be an open subset and $h: U \ra \bR$ a smooth
function which is {\cmssl homogeneous} of degree $k\in \bR^*$, in the sense
that
\be \label{homogEq} \sum_{i=1}^{n+1} x^i\frac{\partial}{\partial x^i}h=kh.\ee
We consider the level set
\[ \mathcal{H} := \{ x\in U \mathbin{|} h(x)=1\} , \]
which we assume nonempty. Notice that if $h$ is not the
zero function we can always rescale $h$ such that this assumption holds.

\bp \label{homogkProp}
For every homogeneous function $h$ as in the above example the inclusion map
\[ \iota: \mathcal{H} \ra \bR^{n+1} \]
is a centroaffine hypersurface embedding with
\[ g = -\frac{1}{k}\iota^*(\tilde\n^2h) , \]
where $(\n , g, \nu)$ are the induced centroaffine data on $\mathcal{H}$.
In particular, $\mathcal{H}\subset \bR^{n+1}$ is nondegenerate if
and only if the Hessian $\tilde\n^2h$ is nondegenerate on $T\mathcal{H}$.
\ep

\pf
By the homogeneity of $h$, $dh(\xi)=k\neq 0$ on $\mathcal{H}$. Therefore
$\mathcal{H}$ is smooth and centroaffine. In order to check the formula for the metric,
let $X$ and $Y$ be vector fields defined on some open subset of $\bR^{n+1}$, which
are tangent to the level sets of $h$. Then on $\mathcal{H}$ we have
\[ \begin{split}
g(X,Y) &= \frac{1}{k}dh (\tilde{\n}_XY)
= \frac{1}{k}\big(\tilde{\n}_X(dhY)-(\tilde{\n}_Xdh)Y\big)\\
&= -\frac{1}{k}(\tilde{\n}_Xdh)Y
= -\frac{1}{k}(\tilde{\n}^2h)(X,Y) .\qedhere
\end{split} \]
\epf

\noindent
Locally every centroaffine hypersurface is defined by a homogeneous function:
\bp \label{homogfProp}
Let $\varphi : M \ra \bR^{n+1}$ be a centroaffine hypersurface immersion, $p\in M$ and $k\in \bR^*$.
Then there exist open neighbourhoods $U'\subset M$ of $p$ and $U\subset \bR^{n+1}$ of $\varphi (p)$
and a smooth homogeneous function of degree $k$ on $U$ such that $\varphi(U') = \{ x\in U \mathbin{|} h(x)=1 \}$.
\ep

\bd \label{centroDef}
A {\cmssl centroaffine structure} on a smooth manifold $M$
is a triple $(\n,g,\nu )$ consisting of a torsion-free connection,
a pseudo-Riemannian metric and a volume form satisfying the following
compatibility conditions:
\begin{enumerate}
\item[(i)] $\n \nu =0$,
\item[(ii)] the curvature tensor $R$ of $\n$ is given by
\[ R(X,Y)Z= -\big(g(Y,Z)X-g(X,Z)Y\big),\quad X,Y,Z\in \mathfrak{X}(M),\]
\item[(iii)] $\n g$ is completely symmetric.
\end{enumerate}
If these conditions are satisfied, $(M,\n , g, \nu )$ is called
a {\cmssl centroaffine manifold}.
The pseudo-Riemannian metric $g$ is called the {\cmssl centroaffine
metric} and the symmetric tensor field $C:=\n g$ is called the {\cmssl cubic
form} of the centroaffine manifold $(M,\n , g, \nu )$.
\ed

\bt  \label{centroThm}\begin{enumerate}
\item[(i)] Let $\varphi : M \ra \bR^{n+1}$ be a nondegenerate centroaffine
hypersurface immersion of a connected manifold $M$
with induced data $(\n , g,\nu )$.
Then $(M,\n , g, \nu )$ is a centroaffine manifold.
\item[(ii)]
Conversely, for a connected and simply connected centroaffine manifold
$(M,\n , g, \nu )$, there exists a centroaffine
hypersurface immersion $\varphi : M \ra \bR^{n+1}$
with induced data $(\n , g,\nu )$.
Furthermore, the immersion $\varphi$ is unique up to
linear unimodular transformations of $\bR^{n+1}$.
\end{enumerate}
\et

\pf
To prove (i) it remains to check
the equations (ii) and (iii) in Definition \ref{centroDef}.
The first equation is obtained by computing the tangent part
of $\tilde{R}(X,Y)Z=0$ with the help of the equations of Gau{\ss} and
Weingarten, where $\tilde{R}$ denotes
the curvature tensor of $\tilde\n$. Similarly, the second equation
is obtained by computing the part proportional to $\xi$.
(These are in fact special cases of the equations of Gau{\ss} and
Codazzi for general hypersurface immersions.)
The statement (ii) can be proven in a similar way as
the fundamental theorem
\cite[Thm.\ 8.1]{NS} of affine differential geometry.
\epf

\subsection{Completeness and closedness of centroaffine hypersurfaces} \label{closeSec}
Our overall approach for proving the completeness of a Riemannian manifold is based on the
following fact, see \cite[Lemma 1]{CHM} for a proof.
\bl
A Riemannian manifold $(M,g)$ is complete if and only if every curve $\gamma : I \rightarrow M$
which is not contained in any compact subset of $M$ has infinite length.
\el

Recall that a submanifold of Euclidean space is called
{\cmssl Euclidean complete} if it is complete with respect to the
Riemannian metric induced by the Euclidean scalar product $\langle \cdot ,\cdot \rangle$.

\bp \label{closedProp}
Let $h : V \ra \bR$ be a smooth homogeneous function of degree $k\neq 0$ defined on some open set $V\subset \bR^{n+1}$
and let $U\subset V$ be an open subset such that $\ol{\mathcal{H}}\subset V$, where
$\mathcal{H} :=  \{ x\in U \mathbin{|} h(x)=1 \}$. Assume that the centroaffine metric of the hypersurface
$\mathcal{H} \subset \bR^{n+1}$  is definite. Then the following hold for every component
$\mathcal{H}_0$ of  $\mathcal{H}$.
\begin{enumerate}
\item[(i)] If $(\mathcal{H}_0,g)$ is complete then $\mathcal{H}_0\subset \bR^{n+1}$ is a closed subset.
\item[(ii)] $\mathcal{H}_0\subset \bR^{n+1}$ is a
closed subset if and only if $\mathcal{H}_0\subset \bR^{n+1}$ is
Euclidean complete.
\end{enumerate}
\ep

\pf
In order to prove (i), let us denote
by $L$ the connected component of the level set $\{ x\in V \mathbin{|} h(x)=1 \}$ which contains $\mathcal{H}_0$.
Thanks to the assumption $\ol{\mathcal{H}}\subset V$, the closure
of $\mathcal{H}_0$ is contained in the level set $\{ x\in V \mathbin{|} h(x)=1 \}$  and thus in $L$. Therefore, if
$\mathcal{H}_0$ is not closed,
then there exists a smooth curve $c : [0,1] \ra L$ such that
$c(0) \in \mathcal{H}_0$ and $c(1)\not\in \mathcal{H}_0$. We can assume without loss of generality that $c([0,1)) \subset
\mathcal{H}_0$. The length of $c$ with respect to the metric $g$ is finite, since it
is given by the integral
\[ \int_0^1 \sqrt{\left|g_{c(t)} (c'(t),c'(t))\right|}dt,\quad g= -\frac{1}{k}\iota^* \tilde{\n}^2h , \]
see  Proposition \ref{homogkProp}, of a continuous function over a compact interval. The continuity of the integrand follows from the
continuity of the Hessian  $\tilde{\n}^2h$ on $V$. This proves (i).

Since the Euclidean length
\[ \int_0^1 \sqrt{\langle c'(t),c'(t)\rangle}dt\]
is also finite, the same argument shows that Euclidean completeness implies that
$\mathcal{H}_0\subset \bR^{n+1}$ is closed. The converse statement follows from the next
simple lemma which finishes
the proof of (ii).
\epf

\bl Let $\varphi : M \ra N$ be an embedding
into a complete Riemannian manifold $(N,g_N)$. If
$\varphi (M) \subset N$ is a closed subset, then
$(M,\varphi^*g_N)$ is complete.
\el

\noindent
{\bf Remark:} The lemma does not extend to injective immersions.

\subsection{Completeness of hyperbolic centroaffine hypersurfaces} \label{auxSec}

In this section we provide some basic results about Euclidean complete connected hyperbolic centroaffine hypersurfaces
$\mathcal{H}\subset \bR^{n+1}$, which will be used in Sections \ref{boundarySec}, \ref{mainSec} and \ref{ApplSec}.  We will first show that the cone
\[ U := \bR^{>0}\cdot \mathcal{H}\]
is open and convex and that  it intersects every affine tangent space $E_p := p + T_p\mathcal{H}$ of $\mathcal{H}$, $p\in \mathcal{H}$,
in a relatively compact convex domain $U\cap E_p \subset E_p$. We will  then parametrize the hypersurface as a radial graph over such a domain and compute the
centroaffine metric in that parametrization as well as its pullback under the radial projection $U\ra \mathcal{H}$. The explicit formulas
involve a positive homogeneous function $h$ on $U$ that is constant on the hypersurface. Building on these preparations, the upshot of this section is Lemma
\ref{concaveLemma}, which provides a criterion for the completeness of the centroaffine metric
in terms of a concavity property of the function $h|_{U\cap E_p}$. All these results require the Euclidean completeness,
with exception of the explicit formulas for the metric, which hold also for local radial parametrizations and projections.
Let us emphasize that Lemma \ref{concaveLemma} will play a key role in the proof of Theorem \ref{mainThm}, which is
the main result of this paper. Another important result  of this section is Proposition \ref {GHProp}, which asserts that $U$ carries a natural
Lorentzian metric that is globally hyperbolic if and only if the Euclidean complete hypersurface $\mathcal{H}$ is complete with respect to the centroaffine metric.

\bp \label{convexProp}
The cone $U$ is open and convex. The map $\phi\colon\bR^{>0}\times\mathcal{H}\to U$ given by $(\lambda,x) \mapsto \lambda x$ is a diffeomorphism.
\ep

\pf
Since $\mathcal{H}$ is centroaffine, $\phi$ is a surjective local diffeomorphism. Thus $U$ is open.

By the Sacksteder--van Heijenoort theorem, see for instance \cite{Wu}, any
Euclidean complete connected hypersurface with positive sectional curvature of dimension $n\ge 2$
is the boundary of  a convex domain $\mathcal{D}$ in $\bR^{n+1}$.
We want to apply this result to our situation.

For this we have to show that
the curvature of $(\mathcal{H},\iota^*\langle.,.\rangle)$ is positive. Note that by Proposition \ref{homogfProp}
we can assume that locally $\mathcal{H}$ is $\{h=1\}$ for some homogeneous function $h\colon U'\to \mathbb{R}$ on some subcone $U'\subset U$.
It follows that (see \cite{GKM}) the second fundamental form of $(\mathcal{H},\iota^*\langle.,.\rangle)$ in
$(\mathbb{R}^{n+1},\langle.,.\rangle)$ is $\frac{1}{|\text{grad }h|}\tilde\nabla^2 h$ with respect to the normal
$-\frac{\text{grad }h}{|\text{grad }h|}$, where $\text{grad }h$ denotes the gradient of $h$ with respect to the Euclidean scalar product $\langle.,.\rangle$.
By our assumptions the second fundamental form is thus definite. Now the Gau{\ss} equation
$$ g_{Euc} (R^{Euc}(X,Y)Y,X)
= \frac{1}{|\text{grad }h|^2}\left(\tilde\nabla^2 h(X,X)\tilde\nabla^2 h(Y,Y)-\tilde\nabla^2 h(X,Y)^2\right)>0$$
for the curvature $R^{Euc}$ of the induced metric $g_{Euc} = \iota^*\langle.,.\rangle$ of the immersion $\iota \colon (\mathcal{H},\iota^*\langle.,.\rangle)$$\to (\mathbb{R}^{n+1},\langle.,.\rangle)$ shows that the sectional curvatures are positive. Hence in the case $n\geq2$, there exists a convex domain $\mathcal{D}\subset\bR^{n+1}$ whose boundary is $\mathcal{H}$.

Notice that $\mathcal{H}$ separates $U$ into two domains
$U^+:=\bR^{>1}\cdot \mathcal{H}$ and $U^-:= \bR^{<1}\cdot \mathcal{H}$ with common
boundary $\mathcal{H}$. The first domain lies on the convex side of $\mathcal{H}$; i.e., every $x\in\mathcal{H}$ has a neighbourhood $W_x$ such that $W_x\cap U^+$ is convex and thus equal to $W_x\cap\mathcal{D}$. This can be seen as follows. The Euclidean second fundamental form of $\mathcal{H}$ with respect to its outer unit normal $\nu$ (with respect to the origin) and the centroaffine metric $g$ on $\mathcal{H}$ differ by a positive conformal factor, because they are affine fundamental forms of $\mathcal{H}$ defined by the Gau{\ss} equation with respect to vector fields $\nu$ and\ $\xi$ that induce the same orientation of $\mathcal{H}$. Since $\mathcal{H}$ is hyperbolic, its Euclidean second fundamental form with respect to $\nu$ is therefore positive definite, which implies the claim that $U^+$ lies on the convex side of $\mathcal{H}$.

Moreover, $U^+$ is contained in $\mathcal{D}$, and $\phi$ is injective. Otherwise some ray $\bR^{>1}\cdot x$ with $x\in\mathcal{H}$ would meet the boundary of $\mathcal{D}$, i.e., it would contain some $x'\in\mathcal{H}$. Since $\mathcal{D}$ is convex, the line from $x$ to $x'$ would lie in $\mathcal{D}$, in contradiction to the fact that the convex side of $\mathcal{H}$ at $x'$ is the outer one.

We claim that $ \mathcal{D}\subset U^+$.
Notice that the origin is not a point of $\mathcal{D}$, because else the convexity of $\mathcal{D}$ and the hyperbolicity of $\mathcal{H}$ would again yield a contradiction.
Thus the line segment connecting any point $p\in \mathcal{D}$ to the origin has to intersect $\mathcal{H}$ in some point $q$.   This implies that
$p\in \bR^{>1}\cdot q\in U^+$, proving that $\mathcal{D}=U^+$. It now follows from the convexity of $\mathcal{D}$ that
$U^+$  and, hence, $U= \bR^{>0}\cdot U^+$ is convex if $n\ge 2$.

In the case $n=1$, $\mathcal{H}\subset (\bR^{\geq1}\cdot x) +T_x\mathcal{H}$ holds for some (in fact, every) $x\in\mathcal{H}$, because otherwise an intermediate value argument (involving $x$ and a hypothetical second element of $\mathcal{H} \cap (x +T_x\mathcal{H})$) would yield a contradiction to the hyperbolicity of $\mathcal{H}$. In particular, $\mathcal{H}$ and, hence, $U$ is contained in the open half-space $\bR^{>0}\cdot x +T_x\mathcal{H} \subset \bR^2$.
This proves the convexity of the connected cone $U$ if $n=1$.

Let $\pi\colon\bR^2\setminus\{0\}\to S^1$ be the projection $x\mapsto\frac{x}{|x|}$ and $\text{incl}_{\mathcal{H}} : \mathcal{H} \ra \bR^2\setminus \{ 0\}$ the inclusion. Since $\phi$ has image contained in a half-space, $\pi\circ\text{incl}_{\mathcal{H}}: \mathcal{H}\to S^1$, which is by centroaffineness an immersion, is not surjective. This implies that $\pi\circ\text{incl}_{\mathcal{H}}$ and thus $\phi$ is injective if $n=1$.

For every $n$, the injectivity of $\phi$ shows now that $\phi$ is a diffeomorphism.
\epf

\noindent
Let us recall that $E_p = p + T_p\mathcal{H}\subset \bR^{n+1}$ denotes the affine hyperplane tangent to $\mathcal{H}$ at $p\in \mathcal{H}$.

\bc The intersection $U\cap E_p\subset E_p$ is a convex domain for all $p\in \mathcal{H}$.
\ec

Next we observe that $\mathcal{H}\subset U$ can be described as the level set of a smooth positive function
$h : U \ra \bR$ homogeneous of degree $k\in \bR^*$. The function $h=h_{k}$ is defined by
\[ h(\l x) := \l^k\quad\mbox{for all}\quad  \l \in \bR^{>0},\; x\in \mathcal{H}. \]
It is well-defined and smooth because $\phi$ considered in Proposition \ref{convexProp} is a diffeomorphism.

We denote by
\[ \psi : U \ra \mathcal{H},\quad x\mapsto \frac{x}{\sqrt[k]{h(x)}}, \]
the radial projection onto $\mathcal{H}$. The restriction
\be \label{varphiEq} \varphi := \psi\big|_{U\cap E_p} : U\cap E_p \ra \mathcal{H}\ee
is a parametrization of the hypersurface. In view of Proposition \ref{homogkProp},
the centroaffine metric of $\mathcal{H}$ in this parametrization is given by
\[ g = -\frac{1}{k}\varphi^*(\tilde{\n}^2h) = -\frac{1}{k}\psi^*(\tilde{\n}^2h)\big|_{U\cap E_p} . \]

\bl \label{psiL}
\[ -\psi^*(\tilde{\n}^2h) = -\frac{1}{h}\tilde{\n}^2h +\frac{k-1}{kh^2}dh^2 . \]
\el

\pf
For $\l >0$ let us denote by $\mu_\l : U \ra U$ the scalar multiplication by $\l$. By homogeneity
\re{homogEq}, we have $\mu_\l^*h=\l^kh$. Since $\mu_\l$ is affine, this implies
\[ \mu_\l^*\tilde{\n}^\ell h = \l^k\tilde{\n}^\ell h \]
for all $\ell \ge 0$. As a consequence, we also have
\[ \tilde{\n}^\ell h\big|_{\l x} = \l^{k-\ell} \tilde{\n}^\ell h\big|_{x} \]
for all  $\ell \ge 0$. In fact,
\[ (\mu_\l^*\tilde{\n}^\ell h)_x (v_1,\ldots , v_\ell )
= \tilde{\n}^\ell h\big|_{\l x}(\l v_1 , \cdots , \l v_\ell )
= \l^\ell \tilde{\n}^\ell h\big|_{\l x}(v_1 , \cdots , v_\ell ) \]
for all  $x\in U$, $v_1,\ldots v_\ell \in T_xU= \bR^{n+1}$.
Next we compute
\[ d\psi_x = h(x)^{-\frac{1}{k}}\mathrm{Id} -\frac{1}{k}h(x)^{-\frac{1}{k}-1}dh_x\ot x.\]
Using these formulas, we calculate
\begin{align*}
(\psi^*\tilde{\n}^2h)_x &= \tilde{\n}^2h_{\psi(x)}(d\psi_x \cdot , d\psi_x \cdot ) = h(x)^{-\frac{k-2}{k}} \tilde{\n}^2h_{x}(d\psi_x \cdot , d\psi_x \cdot )\\
&= h(x)^{-\frac{k-2}{k}}\left[h(x)^{-\frac{2}{k}}\tilde{\n}^2h_x-\frac{2}{k}h(x)^{-\frac{2}{k}-1}dh_x\otimes \tilde{\n}^2h_x(x,\cdot )\right.
\\
&\mspace{102mu}\left.+ \frac{1}{k^2}h(x)^{-\frac{2}{k}-2}\tilde{\n}^2h_x(x,x)dh_x\otimes dh_x\right]\\
&= h(x)^{-1}\cdot \tilde{\n}^2h_x-\frac{2(k-1)}{k} h(x)^{-2} dh_x\otimes dh_x
+\frac{k-1}{k} h(x)^{-3+1} dh_x\otimes dh_x\\
&= h(x)^{-1} \tilde{\n}^2h_x-\frac{k-1}{k} h(x)^{-2} dh_x\otimes dh_x,
\end{align*}
where at the penultimate step we have used that
\be \label{intEqu} \tilde{\n}^2h_x(x,\cdot ) = (k-1)dh_x\ee
in combination with $dh_x(x)=kh(x)$.
The former equation holds because the partial derivatives of first order of $h$ are homogeneous of
degree $k-1$.
\epf

\bc \label{uCor}
The centroaffine metric of the hypersurface $\mathcal{H}$ in the parametrization \re{varphiEq} is given by
\[ g = -\frac{1}{k\bar{h}}\tilde{\n}^2\bar{h} +\frac{k-1}{(k\bar{h})^2}d\bar{h}^2
= -\frac{1}{u}\tilde{\n}^2u , \]
where $\bar{h}$ denotes the restriction of $h$ to $U\cap E_p$, $u:= \sqrt[k]{\bar{h}}$ and $\tilde{\n}$ denotes the flat connection of the affine space $E_p$.
\ec

\noindent
{\bf Remark:} The function $\sqrt[k]{h}= \sqrt[k]{h_k}$ coincides with $h_1$ and is thus independent of $k$.

\bl \label{cpL}
The convex domain $U\cap E_p\subset E_p$ is relatively compact for all $p\in \mathcal{H}$.
\el

\pf
We consider the positive function $u= \sqrt[k]{h}\big|_{U\cap E_p} = h_1\big|_{U\cap E_p}: U\cap E_p \ra \bR$,
which is concave by Corollary \ref{uCor}, since $g$ is positive definite.
Let $B_\d (p)$ be a Euclidean ball in $E_p$, which is relatively compact in
$U\cap E_p$. There exists $\varepsilon >0$ such that
\be \label{estimEq} -\varepsilon \langle \cdot , \cdot \rangle \ge  \tilde{\n}^2u \ee
on $B_\d (p)$.  We will compare $u$ to the concave $C^1$-function $v : E_p \ra \bR$
defined by
\[ v(x) := \begin{cases} 1-\varepsilon |x-p|^2 &\text{if $x\in B_{\delta}(p)$}\\
1+\varepsilon \d^2 -2\varepsilon \d |x-p| &\text{otherwise} . \end{cases} \]
We claim that $v\ge u$ on $E\cap E_p$. We have $v(p)=u(p)=1$
and both functions take their global maximum at $p$, so $dv_p=du_p=0$.
For any point $x\in U\cap E_p\setminus \{p\}$ we consider the line segment $c : [0,1]\ra U\cap E_p$,
$t\mapsto (1-t)p + tx$, from $p$ to $x$. Put
\[ f := (v-u)\circ c:  [0,1]\ra \bR \]
and $t_0 :=  {\min}\big(1,\frac{\d}{|x-p|}\big)$.
We will prove that $f\ge 0$, which implies $v\ge u$, since $x$ was arbitrary.
Notice that $f$ is smooth if $x\in B_\d (p)$. Otherwise $f|_{[0,t_0]}$ and $f|_{[t_0,1]}$
are smooth.
We have $f(0)=f'(0)=0$ and, in virtue of \re{estimEq}, also  $(f|_{[0,t_0]})''\ge 0$.
Using the initial conditions, this implies $f'\ge 0$  on $[0,t_0]$ and, hence, $f\ge 0$ on $[0,t_0]$.
In particular, $f(t_0)\ge 0$ and $f'(t_0)\ge 0$. This proves that $f\ge 0$ if $t_0=1$.
Otherwise it suffices to observe
that $(f|_{[t_0,1]})''=-(u\circ c|_{[t_0,1]})''\ge 0$, which implies
$f'\ge 0$ and finally $f\ge 0$ using the inequalities at $t_0$. So we have
proven that $v\ge u$. As a consequence
\[ U\cap E_p = u^{-1}((0,\infty)) \subset v^{-1}([0,\infty)) \]
and the latter set is compact.
\epf

\bl \label{concaveLemma}
Let $\mathcal{H}\subset \bR^{n+1}$ be a Euclidean complete connected hyperbolic centroaffine hypersurface,
$h : U=\bR^{>0}\cdot \mathcal{H}\ra \bR$ the corresponding homogeneous function of degree $k>0$ and  $p\in
\mathcal{H}$.
Assume that there exists $\varepsilon \in (0,k)$ such that $f=\sqrt[k-\varepsilon]{\bar{h}}$ is concave,
where $\bar{h}= h|_{U\cap E_p}$. Then $\mathcal{H}$ is complete
(with respect to the centroaffine metric $g$).
\el
\pf We first compute
\[ -\frac{1}{k}\tilde{\n}^2f
= \frac{f}{(k-\e)} \left[\left(k-\frac{k}{k-\e}\right)\frac{1}{(k\bar{h})^2}d\bar{h}^2 -\frac{1}{k\bar{h}}\tilde{\n}^2\bar{h}\right] . \]
Comparing with Corollary \ref{uCor}, we see that
\[ g= \frac{k-\varepsilon}{f}\left(  -\frac{1}{k}\tilde{\n}^2f \right) +\frac{\varepsilon}{(k-\varepsilon)(k\bar{h})^2}d\bar{h}^2 . \]
Since the first term is positive semidefinite by assumption, we obtain the estimate
\be \label{estEqu} g \ge \frac{\varepsilon}{(k-\varepsilon)(k\bar{h})^2}d\bar{h}^2 = \frac{\varepsilon}{k^2(k-\varepsilon)} (d\ln \bar{h})^2,
\ee
which implies the completeness as follows. Let $\gamma : I:= [0,T) \ra \mathcal{H}$,  $T\in \bR^{>0} \cup \{\infty\}$, be a curve
which is not contained in any compact subset of $\mathcal{H}$ and $\g_0: I \ra U\cap E_p$ the corresponding curve
in the parametrization $\varphi : U\cap E_p \xrightarrow{\sim} \mathcal{H}$, see \re{varphiEq}. Then there exists a
sequence $t_i \ra T$ such that $\lim_{i\ra \infty}h(\g_0(t_i)) = 0$.
In view of \re{estEqu}, putting $f_0 = h\circ \g_0$, we can estimate the length of $\g$ as follows:
\[ \begin{split}
L(\g)
&\ge \underbrace{\frac{1}{k}\sqrt{\frac{\varepsilon}{k-\varepsilon}}}_{C:=}
\int_0^{t_i} \left|\frac{d}{dt}\ln f_0\right| dt\\
&\ge C\left|\int_0^{t_i} \frac{d}{dt}\ln f_0 \;dt\right|
= C\mathord{\big|\ln f_0(t_i) -\ln f_0(0)\big|} \xrightarrow{i\ra\infty} \infty .\qedhere
\end{split} \]
\epf

\bp \label{GHProp}
Let $\mathcal{H}\subset \bR^{n+1}$ be a Euclidean complete connected hyperbolic centroaffine hypersurface,
$h : U=\bR^{>0}\cdot \mathcal{H}\ra \bR$ the corresponding homogeneous function of degree $k>1$.
Then
\[ g_L := -\frac{1}{k}\tilde{\nabla}^2h \]
is a Lorentzian metric on the convex domain $U$. The Lorentzian manifold $(U,g_L)$ is
globally hyperbolic if and only if $\mathcal{H}$ is
complete (with respect to the centroaffine metric $g$).
\ep

\pf
By the homogeneity of $h$, the position vector $\xi$ satisfies
\begin{align}
g_L(\xi, \xi ) &\stackrel{\phantom{\re{intEqu}}}{=} -(k-1)h<0 \nonumber \\
g_L(\xi, \cdot ) &\stackrel{\re{intEqu}}{=} -\frac{k-1}{k}dh . \label{gradEqu}
\end{align}
The latter equation shows that $\xi$ is perpendicular to the level sets of $h$, on which
the metric $g_L$ restricts to a positive definite metric. Therefore $g_L$ is Lorentzian.

Also due to the homogeneity of $h$, the position vector field $\xi$ on $U$ is a homothetic Killing vector field:
\[ L_\xi g_L = kg_L . \]
The equation \re{gradEqu} shows that it is also a gradient vector field .
Thus
\[ D^L\xi = \frac{k}{2}\mathrm{Id} , \]
where $D^L$ is the Levi-Civita connection of $g_L$.
Rescaling $\zeta :=\frac{2}{k}\xi$ we get  $D^L\zeta = \mathrm{Id}$. Since
the vector field $\zeta$ is obviously complete, this implies that
$(U,g_L)$ is a metric cone:
\[ g_L = -ds^2 +s^2g . \]
Here $U$ is identified with $\bR^{>0}\times \mathcal{H}$ by means of
the diffeomorphism
\[ \bR^{>0}\times \mathcal{H}\ni (s,p)\mapsto \frac{2}{k}sp\in U.\]
With the substitution $s=e^t$ we obtain
\[ g_L = e^{2t}(-dt^2 +g) . \]
This shows that the metric $g_L$ is globally hyperbolic if and only if
the product metric $-dt^2 +g$ on $\bR \times \mathcal{H}$ is.
If $( \mathcal{H},g)$ is complete then the level sets of $t$
are Cauchy hypersurfaces, which implies the global
hyperbolicity by \cite[Cor.\ 39]{O}.  Otherwise there exists an inextendible geodesic
$\g : [0,T) \ra ( \mathcal{H},g)$ of finite length $T$.
This implies that $J(p,q)$ is noncompact if we put $p=(0,\g (0))$, $q=(2T,\g (0))\in \bR \times \mathcal{H}$.
Recall \cite[p.\ 410]{O}  that $J(p,q)$ stands for the smallest set containing all future-pointing causal
curves from $p$ to $q$ and that the sets  $J(p,q)$ are compact
in every globally hyperbolic Lorentzian manifold $M$, for all $p, q\in M$ , see \cite[p.\ 411]{O}.
This shows that $(U,g_L)$ is not globally hyperbolic if $( \mathcal{H},g)$ is incomplete.
\epf

\subsection{Completeness for hyperbolic centroaffine hypersurfaces with regular boundary behaviour} \label{boundarySec}

In this subsection, $\mathcal{H}\subset \bR^{n+1}$ will be always a Euclidean complete connected hyperbolic centroaffine
hypersurface, $U = \bR^{>0}\cdot \mathcal{H}$ and $h : U \ra \bR$ a smooth homogeneous function of degree
$k>1$ such that $\mathcal{H} = \{ p\in U \mathbin{|} h(p)=1 \}$. We assume that $h$ extends to a smooth homogeneous function
$h: V \ra \bR$ defined on some open subset $V\subset \bR^{n+1}$ such that $\ol{U}\setminus\{0\}\subset V$. This assumption is satisfied, for
instance, if the function $h: U \ra \bR$ is polynomial. $0\in \bR^{n+1}$ is excluded in order to keep the level of generality. Note that if a homogeneous
function is smooth at the origin then the degree of homogeneity $k$ is a nonnegative integer. This follows from the fact that all radial derivatives, especially those
of order $n>k$, have to be bounded in $0$, which is not possible for negative or non integer degrees of homogeneity.

\bd  \label{regDef} Under the above assumptions, we say that the hypersurface $\mathcal{H}$ has {\cmssl regular boundary
behaviour} if
\begin{enumerate}
\item[(i)] $dh_p\neq 0$ for all $p\in \partial U\setminus \{ 0\}$. In particular, $\partial U \setminus \{ 0\}$ is smooth.
\item[(ii)] $-\tilde{\n}^2h$ is positive semi-definite on $T(\partial U \setminus \{ 0\})$ with only one-dimensional kernel.
\end{enumerate}
\ed

\noindent
{\bf Example:} The curve $\{ x(x^2-y^2)=1, x>0\}\subset \bR^2$ has regular boundary
behaviour (and is therefore complete by the following theorem),
whereas $\{ x^2y=1, x>0 \} \subset \bR^2$ does not have regular boundary
behaviour (but is still complete). These are the only complete
hyperbolic curves defined by a homogeneous cubic polynomial $h$, up to linear transformations, see \cite[Cor.\ 4]{CHM}.

\bt \label{regThm} Let $\mathcal{H}\subset \bR^{n+1}$ be a Euclidean complete hyperbolic centroaffine
hypersurface with regular boundary behaviour. Then $\mathcal{H}$ is complete (with respect to the
centroaffine metric).
\et

Before we give the proof we would like to discuss the relation of this result to the literature. Melrose \cite[Ch.\ 8]{melrose} considers Riemannian metrics on compact manifolds $M$ with nonempty boundary of the form $x^{2a}dx^2 +x^{2b}H$
where $x\colon M\to [0,\infty)$ is a ``boundary defining function'', i.e. $\{x=0\}=\partial M$ and $dx|_{\partial M}\neq 0$ and $H$
is a smooth tensor field inducing a Riemannian metric on the boundary. Up to a constant factor, the  centroaffine metrics for functions with regular boundary
behaviour correspond to the ``marginally complete'' case $a=-1$, $b=-1/2$ with $x=\bar{h}$ the boundary defining function and $H$ a constant multiple of $-\tilde\nabla^2 \bar{h}$, see Corollary \ref{uCor}.
In this case the map $x\mapsto y=\sqrt{x}$ transforms the metric to a ``conformally compact'' metric $\frac{4}{y^2}\left(dy^2 +\frac{1}{4}H\right)$. The completeness follows from the claim \cite[p.\ 311]{mazzeo} that metrics of this form are complete. For the particular metric considered here this is shown below.

\cite{mazzeo} further claims that each geodesic ray of the conformally compact metric is asymptotic to a single point in the
boundary, the direction approaches the outer normal and the curvature is eventually negative.

\cite{CG} consider metrics as \cite{melrose}, with $a=2b$, but the additional possibility of multiplying the $x^{2a}dx^2$-term by
a smooth function $C$ such that the metric still extends smoothly to the boundary. Under these assumptions and with
$a\le -\frac{1}{2}$ they prove completeness of geodesics whose asymptotic tangents are transversal to the boundary. In \cite{CG2}
the same authors show that $C$ can be assumed constant under certain conditions by an appropriate choice of the boundary function $x$.

As a consequence of Proposition \ref{GHProp}, Theorem \ref{regThm} implies the global hyperbolicity of the Lorentzian metric $g_L$
on $U=\bR^{>0}\cdot \mathcal{H}$ for hypersurfaces $\mathcal{H}$ with regular boundary behaviour, see \cite[Cor.\ 4.34]{FHS} for a related result in
Lorentzian geometry.

\pf
Let $v\in U$ and $E$ a hyperplane trough $v$ which intersects the convex cone $U$ in a relatively
compact domain $B=E\cap U$. Such hyperplanes exist thanks to Lemma \ref{cpL}. We denote
by $\partial B = E\cap \partial U$ the (smooth) boundary of $B$ in $E$. For $\e >0$, we
consider the conical hypersurface $F=F_\e$ which is the union of all the rays
emanating from $\e v$ and intersecting $\partial B$.  It is smooth outside the vertex $\e v$ and so is the
homeomorphism $\psi\big|_{F \cap U} : F \cap U\ra \mathcal{H} $.
Here we recall that $\psi : U \ra \mathcal{H}$ is the map $x \mapsto h(x)^{-1/k}x$.
\bl \label{epsL} There exists $\e_0>0$ such that the  tensor field $-\tilde{\n}^2h$ restricts to a positive definite
metric on a neighbourhood $N_0$ of $\partial B$ in $F_\e $ for all $0<\e < \e_0$.
\el

\pf Let us denote by $\eta\in \mathfrak{X}(V)$ the gradient of $h\in C^\infty (V)$ with respect to the Euclidean scalar product
$\langle \cdot , \cdot \rangle$
in $\bR^{n+1}$. Then, for all $x\in \partial B$,
\[ T_xV=T_x\partial U \oplus \mathrm{span}\{ \eta_x\},\quad  T_x\partial U = T_x\partial B \oplus \mathrm{span}\{ \xi_x\}.\]
The symmetric tensor field $\b = -\tilde{\n}^2h$ on $V$ is positive definite on $T_x\partial B$, by (ii) in Definition  \ref{regDef},
and satisfies
\begin{align*}
\b(\xi_x,\xi_x) &= -k(k-1)h(x) = 0 ,\\
\b(\xi_x,\eta_x) &= -(k-1)dh(\eta_x) = -(k-1)\langle\eta_x,\eta_x\rangle =: c < 0 ,
\end{align*}
and
\[ \b (\xi_x,y)= -(k-1)dh(y)=0,\]
for all $y\in T_x\partial B$, as follows from  the homogeneity of $h$, cf.\ \re{intEqu}.
Now let $(e_3,\ldots ,e_{n+1})$ be a $\beta$-orthonormal basis of $T_x\partial B$,
which we extend by $e_1 := \eta_x$ and $e_2=\xi_x$ to a basis $(e_1,\ldots ,e_{n+1})$ of $T_xV$.
Then
\[ \det (\b (e_i,e_j)_{i,j=1,\ldots, n+1} )= -c^2 \det (\b (e_a,e_b)_{a,b=3,\ldots, n+1}) <0.\]
This shows that $\b_x$ is a Lorentzian scalar product and implies that
the Lorentzian metric $g_L$ is extended by $\frac{1}{k}\b$ to a Lorentzian metric
on a neighbourhood of $\ol{U}\setminus \{ 0\}$ in $V$. Now, to prove the lemma, it suffices to check
for all $x\in \partial B$ that $\b (\nu , \nu )>0$ for a non-zero vector $\nu \in T_xF$ which is orthogonal to the
positive definite hyperplane $T_x\partial B\subset T_xF$ with respect to $\b$.
As such a vector we can take the orthogonal projection of $\e v -x\in T_xF$ onto the
orthogonal complement of $T_x\partial B$ in $T_xF$:
\[ \nu =  \e v-x -\sum_{a=3}^{n+1}\b (\e v -x,e_a)e_a =  \e v-x -\e \sum_{a=3}^{n+1}\b (v,e_a)e_a.\]
Then
\[ \b (\nu , \nu ) = -2\e \b (v,x) + \e^2 \b(v,v)-\e^2 \sum_{a=3}^{n+1}\b (v,e_a)^2.\]
Since $-\b (v,x) = (k-1) dh(v) >0$ and $\partial B$ is compact, we see that $\b (\nu , \nu )>0$ for all
$x\in \partial B$ if $\e$ is sufficiently small.
\epf

Let $\e$, $N_0$ be as in Lemma \ref{epsL}  and put $N=N_0\cap U$, $F=F_\e$.
Then, by Lemma \ref{psiL}, we have
\be \label{hilfEq}(\psi^*g)\big|_N = -\frac{1}{kh}\tilde{\n}^2h\big|_N +\frac{k-1}{(kh)^2}dh^2\big|_N
> \frac{k-1}{(kh)^2}dh^2\big|_N . \ee
This implies the completeness of $g$ as follows. Let $\g : [0, b) \ra \mathcal{H}$, $0<b\le \infty$, be a curve
which is not contained in any compact subset of $\mathcal{H}$ and $\g_F : [0,b)\ra F\cap U$ the corresponding curve in $F \cap U$. Then there exists $\d >0$ and
sequences $t_i \in [0, b)$  and $s_i\in (t_i, b)$ such that $\g_F ([t_i, s_i]) \subset N$,
$h(\g_F(t_i)) >\d$ and $h(\g_F(s_i)) \ra 0$.
In view of \re{hilfEq}, we have
\begin{align*}
L(\g ) &\ge L\big(\g|_{[t_i,s_i]}\big)
=  L(\g_F|_{[t_i,s_i]})
\ge C\int_{t_i}^{s_i}\left| \frac{d}{dt}\ln h \circ \g_F\right| dt\\
&\ge -C\int_{t_i}^{s_i}\frac{d}{dt}\ln h \circ \g_F \;dt
= C \big(\ln h(\g_F(t_i)) - \ln h(\g_F(s_i))\big) \ra \infty ,
\end{align*}
where $C=\frac{\sqrt{k-1}}{k}$.
This shows that $\g$ has infinite length and proves Theorem \ref{regThm}.
\epf

Next we will show that the hypersurfaces with regular boundary behaviour are generic in the
class of hypersurfaces considered in this section. In order to make this statement precise,
let $V \subset \bR^{n+1}$ be an open subset and $k\in (1,\infty)$. We denote
by $\mathcal{F} = \mathcal{F}(V,k)\subset   C^\infty (V)$ the cone consisting of homogeneous functions $h$
of degree $k$ with the property  that there exists an open cone $U\subset V$ such that $\ol{U}\setminus\{0\} \subset V$ and
\[ \mathcal{H}=\mathcal{H}(h,U):= \{ p\in U \mathbin{|} h(p) =1 \} \]
is a Euclidean complete connected hyperbolic centroaffine
hypersurface. (Notice that
for $\mathcal{F}$ to be nonempty $V$ has to contain
an open cone $U$.)
We endow $\mathcal{F}$ with the topology induced by the standard
Fr\'echet topology of $C^\infty (V)$. Recall that the latter is
the coarsest topology for which the semi-norms
$\sup_K \big|\tilde{\nabla}^\ell h\big|$ are continuous for all compact subsets $K\subset V$ and all $\ell = 0,1,\ldots$,
where $|\cdot|$ stands for the Euclidean norm on tensors.

Then we put
\[ \begin{split}
\mathcal{F}_{reg} &\phantom{:}=  \mathcal{F}_{reg}(V,k)\\
&:= \big\{ h\in \mathcal{F} \;\big|\; \text{$\mathcal{H}(h,U)$ has regular boundary behaviour for some $U$ as above} \big\} .
\end{split} \]

\bt \label{genThm} $\mathcal{F}_{reg}$ is a dense open subset of $\mathcal{F}$.
\et

\pf  Let $h\in \mathcal{F}$ and $U \subset V$ an open cone such that $\ol{U}\setminus \{ 0\} \subset V$ and
$\mathcal{H}=\mathcal{H}(h,U)$
is a Euclidean complete connected hyperbolic centroaffine
hypersurface. Replacing $U$ by $\{ p\in U \mathbin{|} h(p)>0 \}$, if necessary, we can assume that  $h>0$ on $U$.
Then $U = \bR^{>0}\cdot \mathcal{H}$. Further let $p\in \mathcal{H}=\mathcal{H}(h)\subset U$ and
$E=E_p$ the affine hyperplane tangent to $\mathcal{H}$ at  $p$. Then we
choose linear coordinates $x_1,\ldots , x_{n+1}$ on $\bR^{n+1}$ such that
$x_1(p)= \cdots =x_n(p)=0$ and $x_{n+1}\big|_E=1$. We claim that
\[ h_\e := h -\e x_{n+1}^k\in \mathcal{F}_{reg}\]
for all $\e\in (0,1)$.

Recall first that $\mathcal{H}$ is closed, by Proposition
\ref{closedProp} (ii). It is mapped to
$\mathcal{H}_\e:= \mathcal{H}(h_\e,U)$ by the following diffeomorphism of the upper half-space
$\{ x_{n+1} >0 \}  \subset \bR^{n+1}$:
\[ x = (x_1,\ldots ,x_{n+1}) = (\vec{x} , x_{n+1})
\mapsto \left(\frac{\vec{x}}{(1+\e x_{n+1}^k)^{1/k}},\; x_{n+1}\right) . \]
This shows that $\mathcal{H}_\e$ is a closed connected smooth hypersurface and therefore Euclidean complete,  again by Proposition
\ref{closedProp} (ii). Next we show that the symmetric tensor field $\b^{\e} := -\tilde{\n}^2h_\e$ on $V$ is
Lorentzian on the cone $U_\e := \bR^{>0}\cdot \mathcal{H}_\e$.  This implies that $\mathcal{H}_\e$
is hyperbolic and, hence, that $h_\e \in \mathcal{F}$.  First we notice that $\b^{\e}$ is Lorentzian when evaluated at
points of $\partial U_\e \setminus \{0\}$.  In fact, since $h_\e-h$ is constant on $E$, $\b^{\e}$ coincides with $\b$ on $E$
and is therefore positive definite on $T\partial B_\e$, where $\partial B_\e$ is the boundary of the
domain $B_\e := U_\e \cap E$ in $E$. Here we are using that $\partial B_\e = \{ p\in E \mathbin{|} h(p)=\e\}$ is
a level set of $h$. On the other hand,  $\xi_x\in T_x\partial U_\e$ is a null vector of $\b^\e$
for all $x\in \partial B_\e$:
\[ \b^\e (\xi_x,\xi_x) = -k(k-1)h_\e (x) = 0 . \]
Moreover, for all
$y\in T_x\partial B_\e$ we have
\[ \b^\e (\xi_x , y) = -(k-1)dh_\e (y) -k(k-1)dh(y) =  0 . \]
Finally, let $\eta^\e$ be the Euclidean gradient of $h_\e$. Then
\[ \b^\e (\xi_x, \eta^\e_x ) = -(k-1)dh_\e (\eta^\e_x) < 0 . \]
As in the proof of Lemma \ref{epsL} , this implies that $\b^\e_x$ has negative determinant and is a Lorentzian scalar product on
$T_xV$ for all $x\in \partial B_\e$.   By homogeneity, the same is true for
all $x\in \partial U_\e \setminus \{ 0\}$. To prove that $\b^\e$ is Lorentzian on $U_\e$ it suffices now to show
that $\det \b^\e$ is negative on $U_\e$.
For all  $x\in B_\e$, we have
\[ \det \b^\e_x = \det \b_x -k(k-1)\e \det A , \]
where the determinant is computed with respect to the basis of $T_xV=\bR^{n+1}$ associated with the coordinates
$x_1,\ldots , x_{n+1}$ and $A$ is the principal $n\times n$-minor obtained by deleting
the last row and column of the Gram matrix of $\b_x^\e$.  Recall that $\b_x$ is Lorentzian and thus  $\det \b_x <0$.
So if $\det A \ge 0$ then it follows that $\det \b^\e_x < 0$ and we are done. Therefore we can assume
that $\det A <0$ and, since $\e < h(x)$ on $B_\e$,
\[ \det \b^\e_x = \det \b_x +k(k-1)\e \left|\det A\right|
< \det \b_x +k(k-1) h(x) \left|\det A\right|
= \det \b_x^{h(x)} . \]
Now observe that $x\in \partial B_{h(x)}$ if $h(x) < \max_Bh = 1$. It follows from the
above discussion that in this case $\b_x^{h(x)}$ is Lorentzian and
$\det \b_x^{h(x)} < 0$. This shows that $\det \b_x^{h(x)} \le 0$ for all
$x\in B_\e$ and implies $\det \b^\e_x <0$ for all $x\in B_\e$. By homogeneity,
this proves that $\b_x^\e$ is a Lorentzian metric on $U_\e$.

Finally, we have to show that $\mathcal{H}_\e$ has regular boundary behaviour.
Since $h_\e = h-\e$ on $E$, we have that
\[ dh\big|_{T_pE} = dh_\e\big|_{T_pE} \quad\text{and}\quad
\b\big|_{T_pE\times T_pE} = \b^\e\big|_{T_pE\times T_pE} \]
for all $p\in V\cap E$ . As $dhT_pE\neq 0$ for all $p\in B\supset B_\e$ and $h_\e$ is homogeneous,
condition (i) in Definition \ref{regDef} is clearly satisfied for all $p\in\partial U_\e\setminus\{ 0\}$. Similarly,
since $\b$ is positive definite on $T_p\partial B_\e = \ker dh|_{T_pE}$
for all $p\in \partial B_\e$ and $h_\e$ is homogeneous, we see that also condition (ii) in Definition \ref{regDef} is satisfied on $T(\partial U_\e \setminus\{0\})$.
\epf

For any integer $k>1$ let us denote by $\mathcal{P}(k) \subset \mathcal{F}(\bR^{n+1},k)$ the subset
consisting of polynomial functions and $\mathcal{P}_{reg}(k) = \mathcal{P}(k)\cap \mathcal{F}_{reg}(\bR^{n+1},k)$.

\bc\label{regcor}
$\mathcal{P}_{reg}(k) \subset \mathcal{P}(k)$ is an open dense subset.
\ec

Next we discuss how many functions with Euclidean complete connected hyperbolic centroaffine level sets there are.

\bt \label{regThm2}
Let $V\subset \bR^{n+1}$ be an open cone. Then for every $k\ge 2$, every compact subset $K\subset V$ and every integer $n\ge 1$ there exists $C=C(k,K,n)<\infty$ such that for all smooth functions
$h\colon V\to\bR$ homogeneous of degree $k$ there exists $h'\in \mathcal{F}_{\text{reg}}(V,k)$ such that
$$\sup_K \big|\tilde\nabla^l h -\tilde\nabla^l h'\big|
\le C \Bigg(\sup_K \sum_{i=0}^{\max (l,2)} \big|\tilde\nabla^i h\big| +1\Bigg)$$
for every $0\le l\le n$.
\et

\pf
We will distinguish two cases, namely $0\notin V$ and $0\in V$.

First assume that $0\notin V$. Let $p\in V$ with $|p|=1$. Consider an open subcone $V'$ of $V$ not containing any nonzero vectors orthogonal to $p$. First
we will construct $h'$ on $V'$ and then extend it to $V$.

Consider an orthonormal basis $\{v_1,\ldots,v_{n+1}\}$ of $\bR^{n+1}$ with $v_{n+1}=p$ and dual basis
$\{\alpha_1,\ldots,\alpha_{n+1}\}$. For $\eta\in\bR$ define the function $h_{\eta}\colon V'\to\bR$, $x\mapsto h(x)+\eta \alpha_{n+1}(x)^k$.
Putting $\eta = |h(p)|+1$ we have $h_\eta(p)>0$.

Next consider the smooth and homogeneous function of degree $k$
$$P\colon\bR^{n+1}\setminus \{0\}\to\bR ,\quad
x\mapsto -\alpha_{n+1}(x)^{k-2}\sum_{i=1}^n\alpha_i(x)^2.$$
Denote with $E_p$ the affine hyperplane in $\bR^{n+1}$ intersecting and orthogonal to $p$, i.e. $E_p=p+\text{span}\{v_1,\ldots,v_n\}$. The second derivative of
$P|_{E_p}^{\phantom{|}}$ at $p$ is negative definite. In fact we have $\tilde{\n}^2\big(P|_{E_p}\big)_p = -2\sum_{i=1}^n \alpha_i\otimes\alpha_i$.

For $\lambda\in\bR$ define $h_{\eta,\lambda} \colon V'\to\bR$, $x\mapsto h_\eta(x)+\lambda P(x)$. Setting $\lambda = \big|\tilde{\n}^2(h|_{E_p})_p^{\phantom{y}}\big| +1$ we conclude that
$\tilde{\n}^2(h_{\eta,\lambda}|_{E_p\cap V})_p^{\phantom{y}} < 0$. Hence $\tilde{\n}^2 (h_{\eta,\lambda}|_{E_p\cap V'})$ is negative definite in a neighbourhood of $p$ in $E_p\cap V'$.
Especially $p$ is a local maximum of $h_{\eta,\lambda}\big|_{E_p}$. For $\lambda$ sufficiently large the set $E_p\cap(h_{\eta,\lambda})^{-1}[0,\infty)$ is contained
in the connected component of $\big\{ q\in V'\cap E_p \mathbin{\big|} \tilde{\n}^2(h_{\eta,\lambda}|_{E_p\cap V'})_q < 0 \big\}$ around $p$. Note that since $\tilde{\n}^2 (h_{\eta,\lambda}|_{E_p\cap V'})_p$
is nondegenerate the map $q\in E_p\cap V'\mapsto d(h_{\eta,\lambda}|_{E_p\cap V'})_q$ is a local diffeomorphism on a neighbourhood $U$ of $p$ in $E_p\cap V'$,
i.e. $d(h_{\eta,\lambda}|_{E_p\cap V'})$ vanishes in $U$ only at $p$. This now implies that the connected component $\mathcal{H}'$ of $h_{\eta,\lambda}^{-1}(1)$ containing
$p$ has regular boundary behaviour. The Euclidean completeness of $\mathcal{H}$ follows from Proposition \ref{closedProp}.

Now consider a smooth function $\mu\colon S^n \to [0,1]$ with support in $V\cap S^n$ and $\mu\big|_{V'\cap S^n}\equiv 1$ such that $\mu \cdot h_{\eta,\lambda}\le 0$ on $(V\setminus V')\cap S^n$ and $\mu(q)=0$ for all $q\in S^n$ orthogonal to $p$.

Define $\hat{h}:\equiv \mu \cdot h_{\eta,\lambda}\big|_{V\cap S^n}$
and consider the canonical $k$-homogeneous extension $h'$ of $\hat{h}$ to $V\setminus \{0\}$. This completes the proof for the case $0\notin V$.

Now assume $0\in V$. Then $k$ is a nonnegative integer. Following the above construction we see that $h-h_{\eta,\lambda}$ is polynomial. In this case we can neglect
the cutoff function $\mu$ and define $h':=h_{\eta,\lambda}$. The claim is now immediate. This completes the proof.
\epf

\bc
For every $k\ge 2$ and every compact set $K\subseteq V$ there exists $C=C(K,k)<\infty$ such that for all polynomials
$h\colon \bR^{n+1}\to\bR$ homogeneous of degree $k$ there exists $h'\in \mathcal{P}_{\text{reg}}(k)$ such that
$$ \sup_K \big|\tilde\nabla^l h-\tilde\nabla^l h'\big| \le C \sup_K \Big(|h|+|\tilde\nabla^2 h|+1\Big) $$
for every $l\ge 0$.
\ec

\section{Projective special real manifolds} \label{2ndSec}

\subsection{Centroaffine structure and intrinsic characterization of projective special real manifolds}\label{intrinsicSec}

\bd \label{PSRDef} A {\cmssl projective special real manifold} is a smooth
hypersurface $\mathcal{H}\subset \bR^{n+1}$ for which
there exists a homogeneous cubic polynomial $h$ on
$\bR^{n+1}$ such that
\begin{enumerate}
\item[(i)] $\mathcal{H}\subset \{h=1\} :=  \{ x\in \bR^{n+1} \mathbin{|} h(x)=1\}$,
\item[(ii)]  the Hessian $\tilde\n^2h$ is negative definite on $T\mathcal{H}$.
\end{enumerate}
\ed

As a consequence of Proposition \ref{homogkProp}, for every projective special real manifold $\mathcal{H}\subset \bR^{n+1}$,
the inclusion $\iota : \mathcal{H}\subset \bR^{n+1}$ is a hyperbolic
centroaffine immersion and, hence, induces a
centroaffine structure $(\n , g ,\nu)$ on $\mathcal{H}$, such that
\be \label{metricEqu} g= -\frac{1}{3}\iota^*(\tilde\n^2h).\ee


\bd \label{intrinsicDef} An {\cmssl intrinsic projective special real manifold} is a
centroaffine manifold $(M,\n,g,\nu)$ with positive definite metric $g$ such that the covariant derivative
of the cubic form $C=\n g$ is given by
\be \label{fundEqu} (\n_X C)(Y,Z,W)= g(X,Y)g(Z,W)+g(X,Z)g(W,Y)+g(X,W)g(Y,Z), \ee
for all $X,Y,Z,W\in \mathfrak{X}(M)$.
\ed

\noindent
{\bf Remark:} The equation \re{fundEqu} implies  that $\n C$ is totally symmetric, that is a
quartic form.

\bt \label{intrinsicThm} \begin{enumerate}
\item[(i)] Let $\mathcal{H}\subset \bR^{n+1}$ be a projective special real
manifold with induced centroaffine structure $(\n , g ,\nu)$. Then
$(\mathcal{H}, \n , g ,\nu)$ is an intrinsic projective special real manifold, that
is satisfies \re{fundEqu}.
\item[(ii)]  Conversely, let $(M,\n,g,\nu)$ be a connected and simply connected
intrinsic projective special real manifold. Then there exists an immersion
$\varphi : M \ra \bR^{n+1}$ such that $\mathcal{H}:= \varphi(M)\subset \bR^{n+1}$ is a
projective special real manifold whose induced centroaffine structure has $\varphi$-pullback $(\n,g,\nu)$.
The immersion $\varphi$ is unique up to linear unimodular transformations of $\bR^{n+1}$.
\end{enumerate}
\et

\noindent
{\bf Remark:} A similar characterization in terms of covariant derivatives of $C$ up to order
$k-2$ can be given for nondegenerate hypersurfaces
which are locally defined by a homogeneous polynomial  $h$ of degree $k\ge 2$.

\pf Let $\mathcal{H}\subset \bR^{n+1}$ be a projective special real
manifold with induced centroaffine structure $(\n , g ,\nu)$.  In order to check \re{fundEqu}, we denote by $H$ the trilinear
form on $\bR^{n+1}$ such that $H(v,v, v)=h(v)$, for all $v\in \bR^{n+1}$.
Differentiating the equation $H(\xi ,\xi ,\xi )=h(\xi)=1$ yields:
\begin{align*}
0 &= H(\xi,\xi,X)\\
0 &= 2H(\xi ,X,Y) + H(\xi , \xi , \n_YX +g(X,Y)\xi ) = 2H(\xi ,X,Y) + g(X,Y)
\end{align*}
for all $X, Y\in \mathfrak{X}(\mathcal{H})$.  Thus
\be g = -2H(\xi ,\cdot ,\cdot )\big|_{T\mathcal{H}\ot T\mathcal{H}}\ee
and, hence,
\[ C(X,Y,Z) = Xg(Y,Z)-g(\n_XY,Z)-g(Y,\n_XZ)=-2H(X,Y,Z) ,\]
for all $X,Y,Z\in \mathfrak{X}(\mathcal{H})$.  Thus
\be C= -2H\big|_{T\mathcal{H}\ot T\mathcal{H}\ot T\mathcal{H}}. \ee
Next we calculate $\n C$ using the previous equations:
\begin{align*}
&(\n_XC)(Y,Z,W)\\
&\stackrel{\phantom{\re{GE}}}{=} XC(Y,Z,W) - C(\n_XY,Z,W)-C(Y,\n_XZ,W)-C(Y,Z,\n_XW)\\
&\stackrel{\phantom{\re{GE}}}{=} -2XH(Y,Z,W)+2H(\n_XY,Z,W)+2H(Y,\n_XZ,W)+2H(Y,Z,\n_XW)\\
&\stackrel{\re{GE}}{=} -2g(X,Y)H(\xi ,Z, W) -2g(X,Z)H(Y,\xi,W)-2g(X,W)H(Y,Z,\xi)\\
&\stackrel{\phantom{\re{GE}}}{=} g(X,Y)g(Z,W)+g(X,Z)g(Y,W)+g(X,W)g(Y,Z),
\end{align*}
for all $X,Y,Z,W\in \mathfrak{X}(\mathcal{H})$. This proves (i).

Let $(M,\n,g,\nu)$ be a connected and simply connected
intrinsic projective special real manifold. Let us denote
by $N=M\times \bR$ the trivial line bundle over $M$,
and by $\xi_0$ its canonical trivializing section.
We claim that \re{fundEqu}
is equivalent to the equation
\[ \tilde\n H =0,\]
where $\tilde\n=\tilde\n^E$ is the flat connection on the vector bundle
$E=TM\oplus N$, which is defined by
\begin{align*}
\tilde\n_XY &:= \n_XY + g(X,Y)\xi_0\\
\tilde\n_X\xi_0 &:= X,
\end{align*}
for all $X,Y\in \mathfrak{X}(M)$ and $H=H_E\in \G (S^3E^*)$ is
defined by
\begin{align*}
H\big|_{TM^{\otimes 3}} &:= -\frac12C=-\frac12\n g\\
H(\xi_0,\cdot, \cdot)\big|_{TM^{\otimes 2}} &:= -\frac12g\\
H(\xi_0,\xi_0, \cdot)\big|_{TM} &:= 0\\
H(\xi_0,\xi_0, \xi_0) &:= 1.
\end{align*}
Let us first show that the curvature $\tilde{R}$ of $\tilde\n$ is zero.
The vanishing of the torsion of $\n$ implies the
equation $\tilde{R}(X,Y)\xi_0=0$
and the equations (ii) and (iii) in Definition \ref{centroDef}
imply $\tilde{R}(X,Y)Z=0$ for all $X,Y,Z\in \mathfrak{X}(M)$.
Next we prove that $\tilde\n H=0$. For $X,Y,Z,W\in \mathfrak{X}(M)$
we compute:
\begin{align*}
(\tilde\n_XH)(\xi_0,\xi_0,\xi_0) &= X\underbrace{H(\xi_0,\xi_0,\xi_0)}_{=1}
-3\underbrace{H(\tilde\n_X\xi_0,\xi_0,\xi_0)}_{=H(X,\xi_0,\xi_0)=0}
= 0\\[1ex]
(\tilde\n_XH)(\xi_0,\xi_0,Y) &= XH(\xi_0,\xi_0,Y) -2H(\xi_0,X,Y) -H(\xi_0,\xi_0,g(X,Y)\xi_0)\\
&= -2H(\xi_0,X,Y)-g(X,Y)
= 0\\[1ex]
(\tilde\n_XH)(\xi_0,Y,Z) &= XH(\xi_0,Y,Z) -H(X,Y,Z) -H(\xi_0,\n_XY,Z) -H(\xi_0,Y,\n_XZ)\\
&= \frac12\Big(-Xg(Y,Z) +C(X,Y,Z) +g(\n_XY,Z) +g(Y,\n_XZ)\Big)
= 0\\[1ex]
(\tilde\n_XH)(Y,Z,W) &= XH(Y,Z,W) -H(\tilde\n_XY,Z,W) -H(Y,\tilde\n_XZ,W) -H(Y,Z,\tilde\n_XW)\\
&= -\frac12(\n_XC)(Y,Z,W) -H(\xi_0,Z,W)g(X,Y) -H(\xi_0,Y,W)g(X,Z)\\
&\mspace{20mu}-H(\xi_0,Y,Z)g(X,W)\\
&= -\frac12\Big( (\n_XC)(Y,Z,W) -g(Z,W)g(X,Y) -g(Y,W)g(X,Z)\\[-0.5ex]
&\mspace{60mu}-g(Y,Z)g(X,W) \Big) .
\end{align*}
This shows that $\tilde\n H=0$ if and only if \re{fundEqu} holds.

Since $(E,\tilde\n )$ is a flat vector bundle over the
simply connected manifold $M$, there exists an isomorphism $\Phi :
(E,\tilde\n)\ra (M\times \bR^{n+1},
\tilde\n)$
identifying $(E,\tilde\n)$ with the trivial vector bundle $(M\times \bR^{n+1},
\tilde\n)$
endowed with its canonical flat connection $\tilde\nabla$.
The restriction $\Phi|_{TM}$ to the subbundle  $TM\subset E$
is a closed vector valued $1$-form  $\phi = (\phi_1, \cdots , \phi_{n+1})$ on $M$.
In fact, for all $X, Y\in \mathfrak{X}(M)$ we have
\[ X\phi (Y)= X\Phi (Y) = \Phi (\tilde\n_XY) =
\phi (\n_XY) + g(X,Y)\Phi (\xi_0)\]
and, hence,
\[ d\phi(X,Y) = X\phi(Y) -Y\phi(X) -\phi([X,Y]) = \phi\big(\n_XY -\n_YX -[X,Y]\big) = 0 .
\]
Since $M$ is simply connected, there exists a smooth map $\varphi : M \ra
\bR^{n+1}$ such that $d\varphi = \phi$. It is a hypersurface immersion because
$\Phi$ being an isomorphism of vector bundles implies that $\phi = \Phi|_{TM}$ is a
monomorphism of vector bundles.  The vector field $\xi := \Phi (\xi_0) : M \ra \bR^{n+1}$
is transversal to $\Phi TM= d\varphi TM$ and verifies
\be \label{WeingEqu} \tilde\n_X\xi = \Phi (\tilde\n_X\xi_0) = \Phi (X) = d\varphi X . \ee
This implies that there exists $v_0\in \bR^{n+1}$ such that
\[ \xi (p) = \varphi (p) +v_0\]
for all $p\in M$.  Composing $\varphi$ with a translation we can
assume that $v_0=0$.   Then $\varphi $ is a centroaffine immersion with
induced data $(\n, g, \det (\xi , \cdots ))$. The induced
volume form $\det (\xi , \cdots )$ is $\n$-parallel (due to \re{WeingEqu}) and,
therefore, coincides with $\nu$ up to a constant factor.
Rescaling $\Phi$, if necessary,  we can assume
that $\nu = \det (\xi , \cdots)$. Now the immersion  $\varphi$ is unique
up to unimodular linear transformation,
by Theorem \ref{centroThm}. Using the identification $\Phi$
of $E$ with the trivial bundle $M\times \bR^{n+1}$,  the parallel section $H=H_E\in \G (S^3E^*)$
corresponds to an element $H=H_{\bR^{n+1}}\in S^3(\bR^{n+1})^*$,
which in turn defines a cubic polynomial $h$ on $\bR^{n+1}$ such that
$h(v)=H(v,v,v)$ for all $v\in \bR^{n+1}$. Now it suffices to show
that $h\circ \varphi =1$, which follows from
\[ 1 = H_E(\xi_0,\xi_0,\xi_0) = H_{\bR^{n+1}} (\xi, \xi, \xi)= h(\xi ) = h(\varphi ).\]
This shows that $\mathcal{H}=\varphi (M)$ is a projective special real manifold.
\epf

\subsection{Relation between completeness and closedness of projective special real manifolds}\label{mainSec}
\bp \label{completeProp} Let $\mathcal{H}\subset \bR^{n+1}$ be a projective
special real manifold with centroaffine metric $g$, see \re{metricEqu}.
Then the following hold.
\begin{enumerate}
\item[(i)] If $(\mathcal{H},g)$ is complete then $\mathcal{H}\subset \bR^{n+1}$ is a
closed subset.
\item[(ii)] $\mathcal{H}\subset \bR^{n+1}$ is a
closed subset if and only if $\mathcal{H}\subset \bR^{n+1}$ is
Euclidean complete.
\end{enumerate}
\ep

\pf Assume that $(\mathcal{H},g)$ is complete or that $\mathcal{H}\subset \bR^{n+1}$ is Euclidean complete.
By taking $V=U=\bR^{n+1}$ in Proposition \ref{closedProp}, we see  that
every component of $\mathcal{H}$ is closed in
$\bR^{n+1}$ and, hence, coincides with one of the finitely many\footnote{The number of
connected components of  a real algebraic set is finite, see \cite{Milnor} and references therein.}  connected components
of $\{ h=1\}$. Then $\mathcal{H}$ is a finite union of closed subsets of
$\bR^{n+1}$ and, therefore, closed. This proves (i) and one of the implications in (ii).
To prove the converse statement in (ii) it is sufficient to remark that the components of the
closed subset  $\mathcal{H}\subset \bR^{n+1}$ are again closed in $\bR^{n+1}$ and, therefore,
Euclidean complete by Proposition \ref{closedProp}.
\epf

\noindent
{\bf Remark:} The previous proposition extends \cite[Prop.\ 5]{CDL}.


\bt \label{mainThm}
Let $\mathcal{H}\subset \bR^{n+1}$ be a projective special real manifold endowed with the centroaffine metric $g$. Then $(\mathcal{H},g)$ is complete if and only if the subset $\mathcal{H}\subset \bR^{n+1}$ is closed.
\et

\pf
In view of Proposition \ref{completeProp} (i), it suffices to show that
a closed projective special real manifold $\mathcal{H}\subset \bR^{n+1}$ is complete.
We can assume without loss of generality that $\mathcal{H}$ is connected, that is a component of the level
set $\{ h=1\}$ of a homogeneous cubic polynomial. By Proposition \ref{completeProp} (ii),
$\mathcal{H}$ is Euclidean complete and is therefore a Euclidean complete connected hyperbolic centroaffine hypersurface
as  considered in Lemma  \ref{concaveLemma}. The unique homogeneous function of degree $k=3$
on $U = \bR^{>0}\cdot \mathcal{H}$ which has the value $1$ on $\mathcal{H}$ coincides with
the restriction of the polynomial $h$ to $U$.  To prove the completeness we will apply Lemma  \ref{concaveLemma}
in the case $k=3$, $\varepsilon =1$. Thus we have to show that the function $\sqrt{h}\big|_{U\cap E}$ is concave,
where $E:=E_p$ is the tangent hyperplane at some point $p\in \mathcal{H}$.   Since $U\cap E \subset E$ is relatively compact (see Lemma \ref{cpL}), for every $x\in U\cap E$ and $v\in T_p\mathcal{H}$ there exists $-\infty <a<b <\infty$ such that the line $x+\bR v\subset E$ intersects the domain $U\cap E$ in the bounded segment
\[ \{ x+tv \mathbin{|} a<t<b \} . \]
We consider the polynomial function $h_0: \bR \ra \bR$ defined by $h_0(t) := h(x+tv)$. It suffices to check that $\sqrt{h_0}''\le 0$ on $(a,b)$. We compute
\[ 4h_0^{3/2}\sqrt{h_0}'' = 2h_0h_0'' -(h_0')^2 \]
and
\[ (2h_0h_0''- (h_0')^2)'= 2h_0h_0'''. \]
Since $h_0'''$ is constant, this shows that the function $f_0 := 2h_0h_0''- (h_0')^2$ is monotone on $(a,b)$. Observing that  $f_0(a) = -(h_0'(a))^2 \le 0$ and $f_0(b) = -(h_0'(b))^2 \le 0$, we see that $f_0\le 0$ on $(a,b)$. This proves that $\sqrt{h_0}'' \le 0$ on $(a,b)$.
\epf

\subsection{Applications}\label{ApplSec}

\bt \label{main_applThm}Let $h$ be a cubic homogeneous polynomial on $\bR^{n+1}$ and $\mathcal{H}$ a locally strictly
convex (i.e.\ definite) component of the level set $\{ h=1\}$. Then $\mathcal{H}\subset \bR^{n+1}$ is a complete
projective special real manifold, which defines a complete
quaternionic K\"ahler manifold of negative scalar curvature diffeomorphic to $\bR^{4n+8}$ by applying first the
$r$-map and then the $c$-map.
\et

\noindent
{\bf Remark:} Notice that the components of the hypersurface
\[ \big\{  x\in \bR^{n+1} \;\big|\; h(x)=1\;\mbox{and}\; g_x\; \mbox{is definite} \big\} \]
are locally strictly convex but are not necessarily
components of the level set $\{ h=1\}$.  In fact, they are in general not closed in the ambient space and therefore incomplete.

\pf
We claim that the centroaffine hypersurface $\mathcal{H}$ is hyperbolic. Assume it is elliptic.

If $n\ge2$ and $\mathcal{H}$ is compact, then Hopf's characterization of ovaloids \cite[p.~122]{H} implies $\bR^{>0}\cdot\mathcal{H} = \bR^{n+1}\setminus\{0\}$. Thus $h$ is positive on $\bR^{n+1}\setminus\{0\}$ and has odd degree, a contradiction.

If $n\ge2$ and $\mathcal{H}$ is noncompact, then the Stoker--Wu theorem \cite{Wu} yields an element $A\in\GL(n+1)$ such that $\mathcal{H}' := A(\mathcal{H}) \subset \bR^n\times\bR$ is the graph of a strictly convex function $f\colon\Omega\to\bR$, where $\Omega$ is a convex open subset of $\bR^n$ and $f$ achieves its minimum at some $x_0\in\Omega$. (To apply the theorem, we used that $\mathcal{H}$ is closed in $\bR^{n+1}$ and thus Euclidean complete by Proposition \ref{closedProp}.) Ellipticity implies that $0\in\bR^{n+1}$ lies in the strict epigraph $C:= \{ (x,y) \mathrel{|} x\in\Omega,\, y>f(x) \}$, which is convex and has boundary $\mathcal{H}'$. The convex set $C\cap(\bR^n\times\{0\})$ is relatively compact in $\bR^{n+1}$, as one sees easily by considering lines from $(x_0,f(x_0))$ to other points on $\mathcal{H}'$, taking the strict convexity of $f$ near $x_0$ and the convexity of $C$ into account. Thus every ray from $0$ in $P:=\bR^n\times\{0\}$ meets $\mathcal{H}'$. We infer that the homogeneous polynomial $h\circ A^{-1}|_{P}$ is positive on $P\setminus\{0\}$ and has odd degree. This is again a contradiction.

If $n=1$, then still $\mathcal{H}$, being a closed embedded centroaffine curve of elliptic type, is the boundary of a convex domain containing the origin in its interior. Therefore, $\mathcal{H} \subset \{h=1\}$ intersects every line through the origin. The set $\{h=0\}$ contains at least one such line, because $h$ has odd degree. Once more, that is a contradiction.

Hence $\mathcal{H}$ is hyperbolic, as claimed. Now the completeness of $\mathcal{H}$ is a consequence of Theorem \ref{mainThm}.
According to \cite[Thm.\ 4]{CHM}, the complete projective special real manifold $\mathcal{H}$ defines a complete projective special K\"ahler domain $M$ of dimension $2n+2$ by the $r$-map. The domain is diffeomorphic to $TU$, where $U= \bR^{>0}\cdot \mathcal{H}$. By Proposition \ref{convexProp}, $U$ is diffeomorphic to a convex domain.
Therefore $M$ is diffeomorphic to $\bR^{2n+2}$. Next, the complete projective special K\"ahler domain $M$ defines a complete quaternionic K\"ahler manifold $N$ of negative scalar curvature by the $c$-map, see \cite[Thm.\ 5]{CHM}. As a differentiable manifold, $N$ is a product $M\times G$, where $G$ is the solvable Iwasawa subgroup of $\SU (1,n+3)$. The latter Lie group is diffeomorphic to $\bR^{2n+6}$.
\epf

\bt
Let $h$ be a cubic homogeneous polynomial on $\bR^{n+1}$ and $\mathcal{H}$ a locally strictly convex component of the level set $\{ h=1\}$. Then
\[ g_L := -\frac{1}{3}\tilde{\nabla}^2h\]
is a globally hyperbolic Lorentzian metric on the convex domain $U= \bR^{>0}\cdot \mathcal{H}$.
\et

\pf
As in the previous theorem, the assumptions imply that  the centroaffine hypersurface $\mathcal{H}\subset
\bR^{n+1}$ is hyperbolic. Thus the result is an immediate consequence of Theorem \ref{mainThm} and Proposition \ref{GHProp}.
\epf

\subsection{An open problem}

For each two natural numbers $n\geq1$ and $k\geq2$, one can consider the following statement:
\begin{description}
\item[$\mathcal{S}(n,k)$:] \textit{For every homogeneous polynomial $h$ of degree $k$ on $\bR^{n+1}$, every locally strictly convex component $\mathcal{H}$ of the level set $\{h=1\}$ is complete with respect to the centroaffine metric.}
\end{description}

Corollary \ref{regcor} and Theorem \ref{regThm} show that for all $n\geq1$ and $k\geq2$, the property in $\mathcal{S}(n,k)$ is true at least for generic polynomials.

\smallskip
As an immediate consequence of Theorem \ref{mainThm}, $\mathcal{S}(n,k)$ is true for all $n\geq1$ in the case $k\in\{2,3\}$:

\bc \label{k23}
Let $h$ be a homogeneous polynomial of degree $k\in\{2,3\}$ on $\bR^{n+1}$ and $\mathcal{H}$ a locally strictly convex component of the level set $\{ h=1\}$. Then $\mathcal{H}\subset \bR^{n+1}$ is complete with respect to the centroaffine metric.
\ec

\pf
The case $k=2$ is trivial since in that case  the tensor field $-\frac12\tilde{\n}^2h$ on $\bR^{n+1}$ inducing the centroaffine metric is constant. The case $k=3$ is part of Theorem \ref{mainThm}.
\epf

Moreover, $\mathcal{S}(1,k)$ is true for every $k\geq2$:
\bt \label{n=1polynom}
Let $h\colon \mathbb{R}^2\to \mathbb{R}$ be a homogeneous polynomial of degree $k\ge 2$ and $\mathcal{H}$ a locally strictly convex connected component of $\{h=1\}$. Then $\mathcal{H}$ is complete with respect to the centroaffine metric.
\et

\noindent
\pf
Since everything is invariant under linear unimodular transformations we can assume that $U=\{x,y>0\}$, i.e. $U$ is the first quadrant in the plane.

Choose a smooth curve $S$ in the first quadrant transversal to the position vector field such that its closure  connects the points $(0,1)$ and $(1,0)$ and is parallel to the $x$-axis near $(0,1)$ and parallel to the $y$-axis near $(1,0)$. Recall from Lemma \ref{psiL} that
$$ -\psi^*(\tilde\nabla^2 h) = -\frac{1}{h}\tilde\nabla^2 h +\frac{k-1}{kh^2}dh\otimes dh $$
on $S\cap U$ for the map $\psi(x,y) = \frac{1}{\sqrt[k]{h(x,y)}}(x,y)$.

For the question of completeness we are only interested in the behaviour of $-\psi^*(\tilde\nabla^2 h)\big|_S$ near $(1,0)$ and $(0,1)$. Again by the invariance under linear unimodular transformations we only need to consider the problem near $(0,1)$. We want to apply the method of Lemma \ref{concaveLemma} for $\varepsilon =1$. Therefore we have to show that $\big(\sqrt[k-1]{h}\big|_S\big)'' \le 0$ near $(0,1)$. Since $S$ is parallel to the $x$-axis in this area the claim follows from
$$ 0 \ge (k-1)h^{\frac{2k-3}{k-1}}\frac{\partial^2}{\partial x^2}{\left(\sqrt[k-1]{h}\right)}
= \frac{2-k}{k-1}\left(\frac{\partial h}{\partial x}\right)^2 +h\frac{\partial^2 h}{\partial x^2} . $$
Note that the right-hand side is polynomial so we only need to consider the monomial $x^l y^{k-l}$ with $l$ minimal appearing in $h$. We know that $1\le l\le k-1$ since $h$ vanishes on both the $x$- and the $y$-axis. Further note that the respective coefficient of $x^ly^{k-l}$ is positive since $h|_U >0$. Then we have
$$ \frac{2-k}{k-1}\big(lx^{l-1}\big)^2 +x^l l(l-1)x^{l-2} = \frac{x^{2l-2}}{k-1}\big(l^2-l(k-1)\big) \le 0 $$
for $x\ge 0$. Now we can use the method of Lemma \ref{concaveLemma} on $S$ near $(0,1)$ and the theorem follows.
\epf

If we consider instead of polynomials the larger class of analytic functions, Theorem \ref{n=1polynom} becomes false, as the following counterexample shows:

\noindent
{\bf Example:} Let $k>1$. The homogeneous function
\[ h(x,y) = \left( \frac{xy}{x+y}\right)^k\]
is real analytic (and even rational if $k$ is an integer) on the quadrant $U :=\{ x>0,y>0\}$.
The hypersurface $\mathcal{H} = \{p\in U \mathbin{|} h(p)=1\}$ is obviously closed in $\bR^2$,
Euclidean complete and can be parametrized
by
\[ \varphi : B \ra \mathcal{H},\quad  p\mapsto h(p)^{-1/k}p,\]
where $B$ is the intersection of
$U$ with the line $\{ x+y=1 \}$.
According to Corollary \ref{uCor}, in this parametrization, the centroaffine metric
is computed from $u = \sqrt[k]{h}\big|_B = xy\big|_B = x(1-x)$ by
\[ g = -\frac{1}{u}\tilde{\n}^2u = \frac{2dx^2}{x(1-x)} . \]
The centroaffine length of the curve $\mathcal{H}$ is thus
\[ \sqrt{2}\int_0^1\frac{dx}{\sqrt{x(1-x)}} < \infty . \]
Since $\mathcal{H}$ is symmetric with respect to the axis $y=x$, this implies that $\mathcal{H}$ is incomplete.

\begin{open}{\ \ }
Given $n\geq2$ and $k\geq4$, decide whether the statement $\mathcal{S}(n,k)$ is true.
\end{open}

Note that it is not possible to generalize the proof of Corollary \ref{k23}, which is based on Lemma \ref{concaveLemma}, to any $k\geq4$. In order to do that, one would have to prove that there exists a constant $c<\frac{k-1}{k}$ such that
\[
c\eta'(x)^2 -\eta(x)\eta''(x) \geq 0
\]
holds for all $x\in[0,1]$ and all polynomials $\eta\colon\bR\to\bR$ of degree $\leq k$ which
\begin{enumerate}[label=(\Alph*)]
\item\label{etaA} satisfy $\eta(0)=\eta(1)=0$ and are positive on the interval $(0,1)$,
\item\label{etaB} satisfy $\frac{k-1}{k}(\eta')^2 -\eta\eta'' > 0$ on $(0,1)$.
\end{enumerate}
Then for $\e = \frac{k}{1-c}\big(\frac{k-1}{k}-c\big)$, the computations in the proof of Lemma \ref{concaveLemma} would show that for each affine line $L$ in $E_p$ which meets $U$, the function $\eta_L:= {\left.h\right|_{L}}$ (which is a polynomial of degree $\leq k$ because $h$ is a polynomial of degree $k$, and which has an affine reparametrization $\eta\colon\bR\to\bR$ with $\eta(0)=\eta(1)=0$ such that $\eta>0$ and $\frac{k-1}{k}(\eta')^2 -\eta\eta'' \geq 0$ hold on $(0,1)$) makes $\sqrt[k-\e]{\eta_L}$ concave on $L\cap U$. This would imply the completeness: $\mathcal{S}(n,k)$ would be true for the considered $k$ and all $n\geq1$.

For $k=3$, the proof of Theorem \ref{mainThm} shows that $c=\frac{1}{2}$ works. For $k=2$, it is easy to see that $c=0$ works. Unfortunately, a constant $c$ with the desired property does not exist for $k\geq4$. The following example demonstrates this in the case $k=4$:

\noindent
{\bf Example:} For each $a\in\bR^{\geq0}$, the fourth-order polynomial
\[
\eta_a(x) := x(1-x)\left(\left(x-\tfrac{3}{20}\right)^2 +\tfrac{51}{20^2} +a\right)
\]
is obviously positive on $(0,1)$ and vanishes at $0$ and $1$. We consider
\[ \begin{split}
P_a(x) := {}&\tfrac{3}{4}\eta_a'(x)^2 -\eta_a(x)\eta_a''(x)\\
= {}&\frac{3\big(14x^2+6x-3\big)^2}{40^2}
+\frac{-80x^4 +188x^3 -42x^2 -24x +9}{40}a
+\frac{4x^2 -4x +3}{4}a^2 .
\end{split} \]
As $14x^2+6x-3$ has precisely one zero in the interval $[0,1]$, namely $x_0:= -\frac{3}{14}+\frac{1}{14}\sqrt{51} \approx 0.2958$, $P_0|_{[0,1]}$ is nonnegative and vanishes precisely at $x_0$.

Since the polynomial $Q:= -80x^4 +188x^3 -42x^2 -24x +9$ is positive at $x_0$, namely $Q(x_0) \approx 2.479$, there exists an $a_0\in\bR^{>0}$ such that for every $a\in(0,a_0]$, $P_a$ is positive on $[0,1]$. (One can even check that $Q$ is positive on $[0,1]$. Since also $4x^2 -4x +3 = 2(x-1)^2 +2x^2 +1$ is positive, $P_a$ is therefore positive on $[0,1]$ for every $a>0$.)

Thus $\eta_a$ satisfies \ref{etaA} and \ref{etaB} with $k=4$ for all $a\in(0,a_0]$. Assume that there exists a constant $c<\frac{3}{4}$ such that $c\eta_a'(x)^2 \geq \eta_a(x)\eta_a''(x)$ holds for all $a\in(0,a_0]$ and $x\in[0,1]$. By continuity, $c\eta_0'(x)^2 \geq \eta_0(x)\eta_0''(x)$ would hold for all $x\in[0,1]$, in particular for $x_0$. We would obtain $c\eta_0'(x_0)^2 \geq \eta_0(x_0)\eta_0''(x_0) = \frac{3}{4}\eta_0'(x_0)^2$, hence $\eta_0'(x_0)=0$. But that is false:
\[
\eta_0'(x_0) = -4x_0^3 +\tfrac{39}{10}x_0^2 -\tfrac{9}{10}x_0 +\tfrac{3}{20} \approx 0.1215 .
\]
Thus for $k=4$, there is no constant $c<\frac{k-1}{k}$ with the desired property described above.


\begin{thebibliography}{GKM}

\newcommand{\book}[1]{\textit{#1}}

\bibitem[CG13]{CG} A.\ \v{C}ap and A.\ R.\ Gover, {\it Projective compactifications and Einstein metrics}, to appear in J. reine angew. Math., DOI 10.1515/crelle-2014-0036, arXiv:1304.1869.

\bibitem[CG14]{CG2} A.\ \v{C}ap and A.\ R.\ Gover, {\it Projective compactness and conformal boundaries}, to appear in Math. Ann., DOI 10.1007/s00208-016-1370-9, arXiv:1406.4225.

\bibitem[CDL]{CDL} V.\ Cort\'es, M.\ Dyckmanns and D.\ Lindemann, {\it Classification of complete projective special real surfaces}, Proc.\ London Math.\ Soc.\ {\bf 109} (2014), no.\ 2, 353--381.

\bibitem[CHM]{CHM} V.\ Cort\'es, X.\ Han and T.\ Mohaupt, {\it Completeness in supergravity constructions}, Comm.\ Math.\ Phys.\ {\bf 311} (2012), no.\ 1, 191--213.

\bibitem[CY]{CY} S.\ Y.\ Cheng and S.-T.\ Yau, {\it Complete affine hypersurfaces. I. The completeness of affine metrics}, Comm.\ Pure Appl.\ Math.\  {\bf 39} (1986), no.\ 6, 839--866.

\bibitem[FHS]{FHS} J.\ L.\ Flores, J.\  Herrera and M.\ S\'anchez, \emph{On the final definition of the causal boundary and its relation with the conformal boundary}, Adv.\ Theor.\ Math.\ Phys.\ {\bf 15} (2011), no.\ 4, 991--1057.

\bibitem[GKM]{GKM} D.\ Gromoll, W.\ Klingenberg and W.\ Meyer, \book{Riemannsche {G}eometrie im {G}ro\ss en}, Lecture Notes in Mathematics 55, Springer-Verlag, Berlin, 1975.

\bibitem[GST]{GST} M.\ G{\"u}naydin, G.\ Sierra  and P.\ K.\ Townsend, \emph{The geometry of $N=2$   Maxwell--Einstein supergravity and Jordan algebras}, Nucl. Phys. {\bf B242} (1984), 244--268.

\bibitem[H]{H} H.\ Hopf, \book{Differential geometry in the large}, Lecture Notes in Mathematics 1000,  Springer-Verlag, Berlin, 1983.

\bibitem[Ma]{mazzeo} R.\ Mazzeo, \emph{The {H}odge cohomology of a conformally compact metric}, J.\ Differential Geometry {\bf 28} (1988), 309--339.

\bibitem[Me]{melrose} R.\ B.\ Melrose, \book{Geometric scattering theory. Stanford Lectures}, Cambridge University Press, Cambridge, 1995.

\bibitem[Mi]{Milnor} J.\ Milnor, \emph{On the Betti numbers of real varieties}, Proc.\ Amer.\ Math.\ Soc.\ {\bf 15} (1964), 275--280.

\bibitem[NS]{NS} K.\ Nomizu and T.\ Sasaki, \book{Affine differential geometry}, Cambridge Tracts in Mathematics 111, Cambridge University Press, Cambridge, 1994.

\bibitem[O]{O} B.\ O'Neill, \book{Semi-Riemannian geometry. With applications to relativity}, Pure and Applied Mathematics, 103, Academic Press, NY, 1983.

\bibitem[W]{Wu} H.\ Wu, {\it The spherical images of convex hypersurfaces}, J.\ Differential Geometry {\bf 9} (1974), 279--290.

\end{thebibliography}
\end{document}